 \documentclass[final,1p,times,twocolumn]{elsarticle}

\usepackage{amssymb}
\usepackage{amstext}
\usepackage{amsmath}
\usepackage{subfigure}
\usepackage[table]{xcolor}
\usepackage{color}
 \usepackage{diagbox}
\usepackage[colorinlistoftodos,textwidth=25mm]{todonotes}
\usepackage{lscape}
 \usepackage{multirow}
 \usepackage{ulem}

\usepackage{amsfonts}
\usepackage{cancel}
\usepackage{rotating}
\usepackage{amssymb}
\usepackage{stackrel}
\usepackage{tikz}
\usetikzlibrary{shapes,snakes}

\newcommand{\Fset}{\mathbb{F}}

\usepackage[amsmath,thmmarks]{ntheorem}
\theoremheaderfont{\bf}
\theorembodyfont{\itshape}
\theoremstyle{plain}
\theoremseparator{:}
\theoremsymbol{}
\newtheorem{lemma}{Lemma}
\newtheorem{theorem}{Theorem}
\newtheorem{corollary}{Corollary}

\newtheorem{remark}{Remark}
\newtheorem{proposition}{Proposition}
\theorembodyfont{\upshape}
\newtheorem{definition}{Definition}
\theoremsymbol{\ensuremath{\Box}}
\newtheorem{example}{Example}

\theoremheaderfont{\sc}
\theorembodyfont{\upshape}
\theoremstyle{nonumberplain}
\theoremsymbol{\ensuremath{\square}}
\newtheorem{proof}{Proof}

 \usepackage[symbol]{footmisc}

\journal{ArXiv}
\begin{document}

\begin{frontmatter}



\title{Constructing Superregular and Block Superregular Matrices}


 \author[label1]{Gustavo Terra Bastos}
 \ead{gtbastos@ufsj.edu.br}
 \address[label1]{ Federal University of São João del-Rei (UFSJ), Brazil\\
 Department of Mathematics and Statistics }

 \author[label2]{Sara D. Cardell \footnote{Corresponding author}}
 \ead{sd.cardell@unesp.br}
 \address[label2]{São Paulo State University (Unesp), Brazil\\
Institute of Geosciences and Exact Sciences (IGCE),
Rio Claro.}

\begin{abstract}
Superregular matrices, i.e., matrices where all square submatrices are non-singular,  have a wide range of applications in communications.
A superregular block matrix  is a broader concept where all full block submatrices, with the appropriate size, are non-singular. 
In this work we propose a construction of block superregular matrices based on the Kronecker product of superregular matrices with non-singular matrices.
Furthermore, we propose two constructions of superregular matrices via  other smaller superregular  matrices over smaller fields.
\end{abstract}


 \end{frontmatter}



\section{Introduction}
Superregular matrices  are also known as MDS matrices and receive such a  name based on   their connection with MDS codes (codes with maximum correction capacity~\cite{Huffman2003bk}).
Besides the theoretical context, these matrices have  applications not only in coding theory, but also in design of modern encryption systems; in particular they are involved in the diffusion layer of block ciphers, for instance, the AES cipher  \cite{Paar2010bk} or the Camellia cipher \cite{Aoki2001}, and in the construction of  hash functions \cite{Stallings2006}.
Some known families of matrices, such as Vandermonde matrices or Cauchy matrices \cite{Roth1989,Roth1985} are known to be superregular, as long as the corresponding field is large enough. 
There are other constructions proposed in the literature, for instance, in \cite{Lacan2003}  the authors introduced a new construction based on the product of matrices with restrictions on the rank. Given an encryption scenario, particular cases of MDS matrices, like involutions and Hadamard matrices, are highly desired for implementation due to the simplicity of computing their inverses and the consequent reduction on implementation costs. 

In general, the construction of superregular matrices over   finite fields is not a easy task due to high computational cost in order to check the determinants of all possible square submatrices.

Superregular block matrices are used in the construction of $\mathbb{F}_q$-linear codes \cite{Louidor2006,Blaum1999,Cardell2013}.
In \cite{Cardell2013}, the authors proposed a construction of block superregular matrices via companion matrices. 
They introduced a ring isomorphism where the image of a superregular matrix is a superregular block matrix composed of powers of the companion matrix of a primitive polynomial.

In this work we use a generalization of that result to do the opposite, that is,  to construct superregular matrices from  superregular block matrices.
Furthermore, we propose a recursive method to obtain superregular block matrices using the Kronecker product. 
 This procedure of constructing superregular matrices over a finite extension allows us to explore new superregular matrices within the same finite field. This technique can be applied to different finite extension fields, proving useful for constructing MDS mappings as a generalization of the subfield construction proposed in~\cite{kamil}. It is noteworthy that there are no constructions of matrices of this type using the Kronecker product, and there are few constructions, not only of superregular matrices but also of block superregular matrices.

Additionally, it is worth mentioning that superregular matrices, whose entries are described by elements of $GL_n (\mathbb{F}_q)$, the group of non-singular matrices over $\mathbb{F}_q$ (referred to as the general linear group), i.e., superregular block matrices, are also essential for encryption procedures. These matrices play a crucial role in constructing lightweight MDS (block) matrices, as discussed in~\cite{zhou}.




This paper is organized as follows. In Section~\ref{sec:prel} we recall the main necessary concepts to understand the rest of the paper. 
In Section~\ref{sec:main}, we introduce a new method to obtain superregular block matrices and superregular matrices using the Kronecker product and an field isomorphism. 
In Section~\ref{sec:var}, we introduce a new method of obtaining superregular matrices using companion matrices of primitive polynomials. 
Finally, he paper concludes in Section~\ref{sec:conc}.

\section{Preliminaries}\label{sec:prel}
In this section we present some basic concepts, starting with the definition of superregular matrix over a finite field.
From now on, we denote by $\mathbb{F}_q$ the Galois field with $q=p^t$ elements, where $p,t \in \mathbb{Z}_{+}$, and $p$ is prime, and $q$ sufficiently large. 

\begin{definition}
A matrix $A$ is said to be a superregular matrix over $\mathbb{F}_q$ if \textbf{every} square submatrix of $A$ is non-singular over $\mathbb{F}_q$.
\end{definition} 

\begin{example}   \label{ex:sup1} 
Consider the matrix $
M=
\left[
\begin{array}{ccc}
   6 &  2  & 2 \\
   4  & 3 & 1 \\
   3  & 3&   4 \\
\end{array}
\right]
$
over $\mathbb{F}_{7}$.
Since all entries, all $2\times 2$    minors (see Table~\ref{tab:ex:sup1}) and the determinant of $M$ are different from zero in $\mathbb{F}_{7}$, we can state that the matrix is superregular over $\mathbb{F}_{7}$.
    \end{example}
  \begin{table}[h!]
      \centering
  \caption{All $2\times 2-$ minors of $M$ in Example~\ref{ex:sup1} \label{tab:ex:sup1}}
\begin{tabular}{|c|c|c|c|}\hline
  \backslashbox{Rows}{Cols}   & $\{1,2\}$ & $\{1,3\}$ & $\{2,3\}$ \\ 
\hline
\multicolumn{1}{|c|}{$\{1,2\}$} &  3 & 5  &3    \\ 
\hline
\multicolumn{1}{|c|}{$\{1,3\}$} & 5  &4  &2   \\ 
\hline
\multicolumn{1}{|c|}{$\{2,3\}$} & 3  & 6  & 2  \\ 
\hline
\end{tabular}
  \end{table}

It is worth mentioning that these matrices are also known as full (or complete) superregular matrices, as well as MDS matrices, due to their relation with MDS (maximum distance separable) codes as we may see in the following theorem.

\begin{theorem}\cite{Huffman2003bk}
An $[n,k,d]$-linear code $C$ over $\mathbb{F}_q$ with generator matrix $G=[I_k,A]$, where $I_k$ is the $k\times k$  identity matrix and $A\in \mathsf{Mat}_{k \times n-k}(\mathbb{F}_{q})$, is MDS if and only if $A$ is a superregular matrix over $\mathbb{F}_q$.
\end{theorem} 

There are several families of matrices that are superregular (as long as the field is large enough), for instance Cauchy and Vandermonde matrices \cite{Roth1985, Burgisser1997bk}. For a encyclopedic reference about superregular/MDS matrix constructions, and their connections with encryption systems, we recommend~\cite{Survey}.

Now, we introduce the definition of block superregular matrix.

\begin{definition}
A matrix $A=\left[A_{ij} \right] \in \mathsf{Mat}_{bm \times bt}(\mathbb{F}_{q})$, where each $A_{ij}$ has size $b\times b$, is said to be a 
$\pmb{b}$-block superregular  matrix if \textbf{every} square submatrix of $A$ consisting of full blocks matrices $A_{ij}$ is nonsigular over $\mathbb{F}_{q}$.
\end{definition}

\begin{example}\label{ex:supblock}
The matrix given by $$  
A=
\left[
\begin{array}{cc|cc}
1&0&1&0\\
0&1&0&1\\ \cline{1-4}
0&1&1&0\\
1&0&1&1\\
\end{array}
\right]
$$
is a 2-block superregular matrix over $\mathbb{F}_2$, since the matrices
$$
\left[
\begin{array}{cc}
    1 &  0\\
   0  & 1
\end{array}
\right],
\left[
\begin{array}{cc}
   0 &  1\\
   1  & 0
\end{array}
\right],
 \left[
\begin{array}{cc}
    1 &  0\\
   1  & 1
\end{array}
\right]
$$
and $A$ itself are non-singular matrices   over $\mathbb{F}_2$. Notice that $A$ is not superregular over $\mathbb{F}_2$.
\end{example}

These matrices are related to MDS $\mathbb{F}_q$-linear codes. 
A code $C$ is said to be  $\Fset_q$-linear of length $n$ over $\Fset_q^b$ if it is a linear code of length $nb$ over $\Fset_q$~\cite{Cardell2013}.
We represent the code $C$ as $C_{\Fset_q^b}$ (respectively $C_{\Fset_q}$) when we consider $C$ as a code over $\Fset_q^b$ (respectively  $\Fset_q$). Both refer to the same set of codewords, we just consider different alphabets in each case. 
\begin{theorem}
[\cite{Blaum1999,Cardell2012PhDT}]\label{th:super:MDS}
Let $
  G
  =[     I_{kb}, A]  
$
be an $kb \times nb$ systematic generator matrix of an $\mathbb{F}_{q}$-linear code $C_{\mathbb{F}_{q}^{b}}$ with parameters $[n,k]$ $\left(\mbox{over }\Fset_q^b \right)$.
Then $C_{\mathbb{F}_{q}^{b}}$ is MDS if and only $A$ is a superregular $b$-block matrix.
\end{theorem}

\begin{example}
Consider the
$\Fset_2$-linear code  $C_{\mathbb{F}_{2}^{2}}$  with parameters $[4,2]$ and generator matrix
\begin{equation*}
G=
[I_4, A]=
\left[
\begin{array}{cc|cc|cc|cc}
1&0&0&0&1&0&1&0\\
0&1&0&0&0&1&0&1\\\cline{5-8}
0&0&1&0&0&1&1&0\\
0&0&0&1&1&0&1&1\\
\end{array}
\right].
\end{equation*}
We saw in Example~\ref{ex:supblock} that $A$ is a superregular 2-block matrix over $\mathbb{F}_2$.
Thus, according to Theorem~\ref{th:super:MDS}, the code  $C_{\mathbb{F}_{2}^{2}}$ is MDS over $\Fset_2^2$.
\end{example}

Let $C\in\mathsf{Mat}_{n}(\mathbb{F}_{q})$ be a companion matrix related to a $n-$degree irreducible polynomial in $\mathbb{F}_q [x]$ (see \cite{Horn1985bk} for more information about companion matrices). Moreover, let $\mathbb{F}_q [C]:=\left\{\sum_{i=0} ^{n-1} f_i C^i : f_i \in \mathbb{F}_q \right\}$, in which $C^0 = I_n.$ Since $\mathbb{F}_{q^n}$ and $\mathbb{F}_q [C]$ are isomorphic representations of the same $n-$degree extension field (see~\cite{Lidl1986bk}), then
$\mathbb{F}_q [C] ^{*}\subset GL_n (\mathbb{F}_q).$

Constructing lightweight superregular/MDS matrices over an $n-$degree extension field, namely, matrices whose $XOR-$count number is small and whose the entries are elements of $GL_n \left(\mathbb{F}_{q} \right)$, is a common research problem on literature - consult, for instance, \cite{zhou,Li2016OnTC}. Hence, from exposed above, it is possible look at constructions of superregular matrices over $GL_n \left(\mathbb{F}_{q} \right)$ as constructions of  superregular $n-$block matrices over $\mathbb{F}_q$.  

At next, a systematization based on the Kronecker product in order to obtain superregular matrices over either $\mathbb{F}_{q^n}$ or $\mathbb{F}_q [C]$ is presented, which allow us to move through both representations freely.

\section{Superregular matrices obtained by Kronecker products}\label{sec:main}

In this section we introduce a method to construct block superregular matrices using the Kronecker product. 
Later on, we used these new matrices to construct superregular matrices in the traditional sense, using a field isomorphism which will be provided by Theorem~\ref{lemmasara}. 

\begin{definition}
Given two matrices  $A\in \mathsf{Mat}_{m \times t}(\mathbb{F}_{q})$  and $B\in \mathsf{Mat}_{n \times s}(\mathbb{F}_{q})$, we define the Kronecker product as: 
\begin{equation*}
A\otimes B=\left[
    \begin{array}{cccc}
      a_{1,1}B& a_{1,2}B&\ldots & a_{1,t}B  \\
      a_{2,1}B& a_{2,2}B&\ldots & a_{2,t}B  \\
      \vdots &\vdots &&\vdots\\
            a_{m,1}B& a_{m,2}B&\ldots & a_{m,t}B  \\   \end{array}\right].
\end{equation*}

If $m=t$ and $n=s$, the determinant of $A\otimes B$ is given by (see, for example, \cite{Horn1994bk}):
    \begin{equation}\label{eq:det:kron}
    \det(A\otimes B)=(\det(A))^n(\det(B))^m .
    \end{equation}
 \end{definition}
 

\begin{lemma}[Schur complement]
Let   $A$, $B$, $C$ and $D$ be matrices of sizes $p\times p$, $p\times (m-k)$, $(m-k)\times p$ and $(m-p)\times (m-p)$, respectively.
If $\det(A)\not = 0$ then 
\[
\det
\left(
\begin{array}{cc}
A & B\\
C & D\\
\end{array}
\right)
=\det(A)\det(D-CA^{-1}B)
\]
\end{lemma}

\begin{example}\label{ex:1}
Consider
$
A=\left[
\begin{array}{cc}
  1    & 2  \\
3 & 4 
\end{array}
\right]
$
superregular over $\mathbb{F}_7$ and
$
B=\left[
\begin{array}{cc}
  1    & 1  \\
   0 & 3 
\end{array}
\right].
$
Since $B$ is non-singular, the matrix 
$$
M=A\otimes B=\left[
\begin{array}{cc|cc}
    1 &  1 &  2  & 2\\
   0  & 3  & 0 &  6\\\hline
   3  & 3 &  4 &  4\\
   0  & 2  & 0 &  5\\
\end{array}
\right]
$$
is a 2-block superregular matrix over $\mathbb{F}_7$ and the determinant is given by
$$\det(M)=\det(A)^2\det(B)^2=5^2\cdot 3^3\equiv 1\bmod 7.$$
Notice that $A\otimes B$ is not superregular. 
\end{example}

At this point it is natural to ask what happens if, for instance, $A$ is a non-square matrix. In that case, we cannot use the Schur complement anymore. However, there are other tools we can use it to generalize it.

\begin{theorem}\label{thm:1}
Consider $A\in \mathsf{Mat}_{m \times t}(\mathbb{F}_{q})$ a superregular matrix and $B\in \mathsf{Mat}_{n \times n}(\mathbb{F}_{q})$ a non-singular  matrix.
The matrix $M=A\otimes B \in \mathsf{Mat}_{mn \times tn}(\mathbb{F}_{q})$ is an $n$-block superregular matrix. 
\end{theorem}

\begin{proof}
Given $M=A\otimes B \in \mathsf{Mat}_{mn \times tn}(\mathbb{F}_{q})$, take the random square submatrix $S$ 
\begin{equation*}
S=\left[\begin{array}{cccc}
    a_{i_1,j_1}B &a_{i_1,j_2}B &\ldots &a_{i_1,j_b}B \\
     a_{i_2,j_1}B &a_{i_2,j_2}B &\ldots &a_{i_2,j_b}B \\
     \vdots &\vdots &&\vdots \\
     a_{i_b,j_1}B &a_{i_b,j_2}B &\ldots &a_{i_b,j_b}B 
\end{array}\right]_{bn \times bn}    
\end{equation*}
of $M$, where $b\leq \min\{m,t\}$. 
Notice that $S=A_{i_b , j_b}\otimes B$, where $$A_{i_b , j_b}= \left[\begin{array}{cccc}
    a_{i_1,j_1} &a_{i_1,j_2} &\ldots &a_{i_1,j_b} \\
     a_{i_2,j_1} &a_{i_2,j_2} &\ldots &a_{i_2,j_b} \\
     \vdots &\vdots &&\vdots \\
     a_{i_b,j_1} &a_{i_b,j_2} &\ldots &a_{i_b,j_b}
\end{array}\right]$$ is a submatrix of $A$.
Then the result follows once $A$ is  superregular, $B$ is non-singular, and (according to expression~(\ref{eq:det:kron})) $\det(S)=(\det(B))^b (det\left(A_{i_b , j_b}\right))^n\neq 0$ .
\end{proof}

In general, the Kronecker product of superregular matrices is not a superregular one. The following result extends the Theorem~\ref{thm:1}.

\begin{theorem}\label{teokronprod}
Let $A_i\in \mathsf{Mat}_{n_i \times n_i}(\mathbb{F}_{q})$ be square superregular matrices, for $1\leq i \leq l $, 
and $B\in \mathsf{Mat}_{n \times n}(\mathbb{F}_{q})$, a non-singular square matrix. Then, the matrix 
\begin{equation*}
M =  A_1\otimes A_2 \otimes \stackrel{l}{\cdots }\otimes A_l \otimes B 
\end{equation*}
is a $\frac{N  n}{n_1}$-superregular matrix of size $N  n\times N  n$ over $\Fset_q$, where $\displaystyle{N=\prod_{i=1}^{l}n_i}$.
\end{theorem}

\begin{proof}
According to Theorem~\ref{thm:1}, $M_1=A_l\otimes B$ is a $n$-block superregular matrix of size $n_ln\times n_ln $.
Now, we know that $M_1$ is non-singular, by definition of block superregular. 
Therefore, we can apply  Theorem~\ref{thm:1} again and state that $M_2=A_{l-1} \otimes  M_1$ is $n_ln$-superregular  of size $n_{l-1}n_{l}n\times n_{l-1}n_{l}n$.
We can apply this argument $l$ times and obtain that
$M_l =A_1\otimes A_2\otimes \stackrel{l}{\cdots}\otimes A_l\otimes B$ is $\frac{Nn}{n_1}$-block superregular matrix of size $Nn\times Nn$.
\end{proof}

\begin{corollary}
\label{thm:3}
Consider $A_i\in \mathsf{Mat}_{n \times n}(\mathbb{F}_{q})$   superregular matrices for $1\leq i \leq l $ and $B\in \mathsf{Mat}_{n \times n}(\mathbb{F}_{q})$ a non-singular matrix. The matrix $M_l =A_1\otimes A_2\otimes \stackrel{l}{\cdots}\otimes A_l\otimes B$ is a superregular $n^{l}$-block matrix of size $n^{l+1}\times n^{l+1}$ over $\Fset_q$. 
\end{corollary}

\begin{example}
Consider again the matrices $A$ and $B$
from Example~\ref{ex:1}.
We can construct now the matrix:
$$
M_2=A\otimes M_1=A\otimes A\otimes B=
\left[ 
\begin{array}{cccc|cccc}
   \textbf{1}&\textbf{1}&\textbf{2}&\textbf{2}&2&2&4&4\\
   \textbf{0}&\textbf{3}&\textbf{0}&\textbf{6}&0&6&0&5\\
   3&3&4&4&6&6&1&1\\
   0&2&0&5&0&4&0&3\\\hline
   \textbf{3}&\textbf{3}&\textbf{6}&\textbf{6}&4&4&1&1\\
   \textbf{0}&\textbf{2}&\textbf{0}&\textbf{4}&0&5&0&3\\
   2&2&5&5&5&5&2&2\\
   0&6&0&1&0&1&0&6\\
\end{array}
\right]
$$
It is possible to check that  $M_2$ is a $4$-block superregular matrix in $\Fset_7$.
Notice that $M_2$ is not $2$-block superregular, since for example, the full block submatrix $
\left[ 
\begin{array}{cc|cc}
  1&1&2&2\\
  0&3&0&6\\\hline
  3&3&6&6\\
  0&2&0&4\\
\end{array}
\right] 
$
is singular. 
\end{example}

\begin{theorem}\label{th:B(B)}
If $A\in \mathsf{Mat}_{ r\times s}(\mathbb{F}_{q})$ is superregular and $B, B_1, B_2,\ldots, B_s \in \mathsf{Mat}_{ n}(\mathbb{F}_{q})$ are non-singular, then 
\begin{equation}\label{mat_prod_b}
A\otimes B\left(B_1 , B_2 , \ldots ,B_s \right)=\left[\begin{array}{cccc}
    a_{11}B\cdot B_1 &a_{12}B\cdot B_2 &\ldots &a_{1s}B\cdot B_s \\
     a_{21}B\cdot B_1 &a_{22}B\cdot B_2 &\ldots &a_{2s}B\cdot B_s \\
     \vdots &\vdots &\ddots&\vdots \\
     a_{r1}B\cdot B_1 &a_{r2}B\cdot B_2 &\ldots &a_{rs}B\cdot B_s \\
\end{array}\right]_{rn \times sn}   
\end{equation}
is $n$-block superregular.
\end{theorem}

\begin{proof}
Let
\begin{equation*}
S=\left[\begin{array}{cccc}
    a_{i_1,j_1}B \cdot B_{j_1}&a_{i_1,j_2}B\cdot B_{j_2} &\ldots &a_{i_1,j_b}B \cdot B_{j_b}\\
     a_{i_2,j_1}B \cdot B_{j_1} &a_{i_2,j_2}B\cdot B_{j_2} &\ldots &a_{i_2,j_b}B \cdot B_{j_b}\\
     \vdots &\vdots &&\vdots \\
     a_{i_b,j_1}B\cdot B_{j_1} &a_{i_b,j_2}B\cdot B_{j_2} &\ldots &a_{i_b,j_b}B\cdot B_{j_b} 
\end{array}\right]_{bn \times bn}    
\end{equation*}
be a submatrix of $A\otimes B\left(B_1 , B_2 , \ldots, B_s\right)$, where $\left\{i_1 , i_2 , \ldots, i_b \right\}\subseteq \{1,2,\ldots,r\}$, and $\left\{j_1 , j_2 , \ldots, j_b \right\}\subseteq \{1,2,\ldots,s\}$. Since
\begin{equation*}
S=\left(A_{i_b , j_b}\otimes B\right)\cdot diag\left(B_{j_1}, B_{j_2},\ldots,B_{j_b}\right), 
\end{equation*}
for $A_{i_b , j_b}= \left[\begin{array}{cccc}
    a_{i_1,j_1} &a_{i_1,j_2} &\ldots &a_{i_1,j_b} \\
     a_{i_2,j_1} &a_{i_2,j_2} &\ldots &a_{i_2,j_b} \\
     \vdots &\vdots &&\vdots \\
     a_{i_b,j_1} &a_{i_b,j_2} &\ldots &a_{i_b,j_b}
\end{array}\right]$ a submatrix of $A$, then the result follows from Theorem~\ref{thm:1}, once 
\begin{equation*}
\det(S)=\det\left(A_{i_b , j_b}\otimes B\right)\cdot det\left(B_{j_1}\right)\cdot det\left(B_{j_2}\right)\cdot \ldots \cdot det\left(B_{j_b}\right)\neq 0
\end{equation*}
\end{proof}

\begin{corollary}\label{co:Bi}
    If $A\in \mathsf{Mat}_{ r\times s}(\mathbb{F}_{q})$ is superregular and $ B_1, B_2,\ldots, B_s \in \mathsf{Mat}_{ n}(\mathbb{F}_{q})$ are non-singular, then the matrix  given by
\begin{equation}\label{mat_prod_b2}
M= \left[\begin{array}{cccc}
    a_{11} \cdot B_1 &a_{12} \cdot B_2 &\ldots &a_{1s} \cdot B_s \\
     a_{21} \cdot B_1 &a_{22} \cdot B_2 &\ldots &a_{2s} \cdot B_s \\
     \vdots &\vdots &\ddots&\vdots \\
     a_{r1} \cdot B_1 &a_{r2} \cdot B_2 &\ldots &a_{rs} \cdot B_s \\
\end{array}\right]_{rn \times sn}   
\end{equation}
is $n$-block superregular.
\end{corollary}
\begin{proof}
    The result follows from Theorem~\ref{th:B(B)}
    taking $B=I$.
\end{proof}

Next result - Theorem~\ref{thm:novamatriz} - is a direct consequence of Theorem~\ref{teokronprod} and provides a construction of superregular matrices via block superregular matrices and companion matrices.
First, we need to introduce the next Theorem. 

\begin{theorem}\cite{Cardell2012PhDT}\label{lemmasara}
Let $\alpha \in \mathbb{F}_{q^n}$ be a fixed primitive element, and $C\in\mathsf{Mat}_{ n}(\mathbb{F}_{q})$ the companion matrix of an $n$-degree  primitive polynomial in $\mathbb{F}_q [x]$.
The map $\psi: \mathbb{F}_{q^n} \rightarrow \langle C \rangle=\{O, I_n, C, C^2, \ldots, C^{q^n-2} \}$,
 with $\psi(\alpha)=C$, is a field isomorphism. 
 Furthermore, it is possible to extend this isomorphism in order to obtain the ring isomorphism  
\begin{equation}\label{inversemapping}
\begin{array}{cccc}
\Psi: &\mathsf{Mat}_{t\times n}(\mathbb{F}_{q^n}) &\rightarrow  &\mathsf{Mat}_{t\times n}(\langle C \rangle)  \\
&A=[a_{ij}]&\mapsto& \Psi(A)=\left[\psi(a_{ij})\right].
\end{array}    
\end{equation}
\end{theorem}

\begin{theorem}\label{thm:novamatriz}
Let $A_i\in \mathsf{Mat}_{n_i \times n_i}(\mathbb{F}_{q})$ be  superregular matrices, for $1\leq i \leq l $ and $C\in \mathsf{Mat}_{n \times n}(\mathbb{F}_{q})$, 
the companion matrix of an $n$-degree primitive polynomial $p(x) \in \mathbb{F}_q [x]$.
Then, the matrix 
\begin{equation*}
M =  A_1\otimes A_2 \otimes \stackrel{l}{\cdots }\otimes A_l \otimes  C 
\end{equation*}
is a $\frac{N n}{n_1}$-block superregular matrix of size $N  n\times N  n$ over $\Fset_q$, where $\displaystyle{N=\prod_{i=1}^{l}n_i}$. Furthermore, the matrix $\Psi^{-1}(M)$ is a new  $N\times N$ superregular matrix over $\mathbb{F}_{q^n}$.
\end{theorem} 

\begin{proof}
The first part follows directly from Theorem~\ref{teokronprod}. Now, given the $\frac{Nn}{n_1}$-block superregular matrix $\displaystyle{M =  A_1\otimes A_2 \otimes \stackrel{l}{\cdots }\otimes A_l \otimes  C   \in \mathsf{Mat}_{Nn \times Nn}(\mathbb{F}_{q})}$, notice that such a matrix might be represented by:

\begin{equation*}
M = \left[\begin{array}{cccc}
     \alpha_{11} C & \alpha_{12} C  &\ldots &\alpha_{1N} C  \\ 
    \alpha_{21} C  & \alpha_{22} C  & \ldots &\alpha_{2N} C  \\
      \vdots& \vdots& &\vdots\\
    \alpha_{N1} C  &\alpha_{N2} C  &\ldots &\alpha_{NN} C 
\end{array}\right]_{N \times N}=\left[\begin{array}{cccc}
      C^{i_{11}}& C^{i_{12}} &\ldots &C^{i_{1N}} \\ 
    C^{i_{21}} & C^{i_{22}} & \ldots &C^{i_{2N}} \\
   \vdots& \vdots& &\vdots\\
    C^{i_{N1}} &C^{i_{N2}} &\ldots &C^{i_{NN}}
\end{array}\right]_{N \times N} 
\end{equation*}
From the map $\Psi^{-1}$, it is straightforward  that
\begin{equation*}
\Psi^{-1}(M) = \left[\begin{array}{cccc}
      \alpha^{i_{11}}& \alpha^{i_{12}} &\ldots &\alpha^{i_{1N}} \\ 
    \alpha^{i_{21}} & \alpha^{i_{22}} & \ldots &\alpha^{i_{2N}} \\
   \vdots  & \vdots& &\vdots\\
    \alpha^{i_{N1}} &\alpha^{i_{N2}} &\ldots &\alpha^{i_{NN}}
\end{array}\right]_{N \times N},
\end{equation*}
is superregular over $\mathbb{F}_{q^n}$.
\end{proof}

\begin{remark}
Notice that the matrix $C$ has order $q^n-1$, and each matrix $C^t$, where $\gcd(t,q^n-1)=1$, can be also used to construct a field isomorphism similar to $\Psi$.
\end{remark}

\begin{example}\label{ex:sup}
Let $A = \left[\begin{array}{ccc}
     1 & 2 & 2 \\
     2 & 1 & 3\\
     3 & 2 & 4
\end{array}\right]\in \mathsf{Mat}_{ 3}(\mathbb{F}_{5})$ be a superregular matrix over $\Fset_5$ and 
$C=\left[\begin{array}{ccc}
0 & 0 & 2\\
1 & 0 & 2\\
0 & 1 & 0
\end{array}
\right]\in \mathsf{Mat}_{3}(\mathbb{F}_{5})$ the companion matrix of the primitive polynomial $p(x)=x^3 +3x +3$.
Consider the matrix 

\begin{equation*}
A\otimes C=\left[\begin{array}{ccc|ccc|ccc}
0 & 0 & 2 & 0 & 0& 4 &0 &0& 4\\ 
1 & 0 & 2 & 2 & 0& 4 &2 &0& 4\\
0 & 1 & 0 & 0 & 2& 0 & 0 & 2 & 0\\\hline
0 & 0 & 4  & 0 & 0& 2 &0 &0& 1\\ 
2 & 0 & 4 & 1 & 0& 2 &3 &0& 1\\ 
0 & 2 & 0 & 0 & 1 & 0 & 0 &3& 0\\\hline
0 & 0 & 1 & 0 & 0& 4 &0 &0& 3\\
3 & 0 & 1 & 2 & 0 & 4 &4 &0& 3\\
0 & 3 & 0 & 0 & 2& 0 &0 &4& 0 
\end{array}\right]\in \mathsf{Mat}_{ 9}(\mathbb{F}_{5}).
\end{equation*}

Since $p(x)$ is a primitive polynomial over $\mathbb{F}_{5}$, the matrix $C$ is non-singular over $\Fset_5$ and has order $5^3-1$.
Furthermore, the blocks of $A\otimes C$ can be considered as powers of $C$ and
\begin{equation*}
M= A \otimes C =  \left[\begin{array}{cccc}
    C & C^{32} & C^{32}\\
    C^{32} & C  & C^{94}  \\
    C^{94} & C^{32} & C^{63} 
\end{array}\right] \in \mathsf{Mat}_{3}(\langle C\rangle)
\end{equation*}
is a $3-$block superregular over $\mathbb{F}_5$. Therefore, applying $\Psi^{-1}$ (the inverse of (\ref{inversemapping})) to $M$, we have

\begin{equation*}
\Psi^{-1}(M)=\left[\begin{array}{ccc}
    \alpha & \alpha^{32} & \alpha^{32} \\
    \alpha^{32} & \alpha & \alpha^{94} \\
    \alpha^{94} & \alpha^{32} & \alpha^{63} 
\end{array}\right]    \in \mathsf{Mat}_{3}(\mathbb{F}_{5^3}) 
\end{equation*}
where $\alpha$ is a primitive element of such finite extension field. According to  Theorem~\ref{thm:novamatriz}, such matrix is superregular over $\mathbb{F}_{5^3}$.

\begin{table}[h!]\label{tabela}
\centering
\caption{All $2\times 2-$ minors of $\Psi^{-1}\left(M\right)$ }
\label{tab:ex:sup}
\begin{tabular}{|c|c|c|c|}\hline
  \backslashbox{Rows}{Cols}   & $\{1,2\}$ & $\{1,3\}$ & $\{2,3\}$ \\ 
\hline
\multicolumn{1}{|c|}{$\{1,2\}$} &$2\alpha^2$  & $4\alpha^2$ & $4\alpha^2$  \\ 
\hline
\multicolumn{1}{|c|}{$\{1,3\}$} & $\alpha^2 $ & $3\alpha^2$ & $4\alpha^2$\\ 
\hline
\multicolumn{1}{|c|}{$\{2,3\}$} & $\alpha^2$  & $4\alpha^2$ & $3\alpha^2$ \\ 
\hline
\end{tabular}
\end{table}

Notice that the matrix $\Psi^{-1}(M)$ is a multiple of $A$ over $\mathbb{F}_{5^3}$ (in fact $\Psi^{-1}(M)=\alpha A$). The determinant of $\Psi^{-1}(M)$ is $\det(\Psi^{-1}(M))=\alpha^3 \det(A)$, and all determinants of the possible $2\times 2-$submatrices over $\mathbb{F}_{5^3}$ are described at Table~\ref{tab:ex:sup}.
\end{example}

From the different representations of $\mathbb{F}_{q^n}$, along with the matrix isomorphism given by Theorem~\ref{lemmasara}, allow us to provide other superregular matrices over $\Fset_{q^n}$.

\begin{example}\label{ex22}
Let $q(x)=x^3 + 3x +2 \in \mathbb{F}_5 [x]$ be a primitive polynomial over $\mathbb{F}_5$ and 
 $$D =\left[\begin{array}{ccc}
    0 & 0 & 3 \\
    1 & 0 & 2 \\
    0 & 1 & 0 
\end{array}\right] \in \mathsf{Mat}_{ 3}(\mathbb{F}_{5})$$
its corresponding companion matrix.
Let $\gamma \in \mathbb{F}_{5^3}$ be a primitive element which is a root of $q(x)$. 
Consider also  the primitive polynomial given in the Example~\ref{ex:sup}, $p(x)=x^3 +3x +3$, and its corresponding companion matrix $C$. From the uniqueness (up to isomorphism) of finite fields, it is natural to establish the following correspondences $\alpha \mapsto \gamma$ and $C \mapsto D$, and extend them to (by abuse of notation) the respective ring isomorphism~\eqref{inversemapping}.
So, we have
\begin{equation*}
N = A \otimes D =\left[\begin{array}{ccc}
    D & D^{94} & D^{94} \\
    D^{94} & D & D^{32} \\
    D^{32} & D^{94} & D^{63} 
\end{array}\right] \mbox{ and }\Psi^{-1}\left(N\right) =  \left[\begin{array}{ccc}
    \gamma & \gamma^{94} & \gamma^{94} \\
    \gamma^{94} & \gamma & \gamma^{32} \\
    \gamma^{32} & \gamma^{94} & \gamma^{63} 
\end{array}\right].
\end{equation*}
Consequently, from the correspondence $\alpha \mapsto \gamma,$ one may state that $ \left[\begin{array}{ccc}
    \alpha & \alpha^{94} & \alpha^{94} \\
    \alpha^{94} & \alpha & \alpha^{32} \\
    \alpha^{32} & \alpha^{94} & \alpha^{63} 
\end{array}\right] \in \mathsf{Mat}_{ 3}(\mathbb{F}_{5})$
is also a superregular matrix. 
\end{example}

Recall that all distinct $\mathbb{F}_q -$automorphisms $\sigma_i$ of $\mathbb{F}_{q^n}$ are exactly the mappings $ \sigma_i \left(x\right) = x^{q^i}$, for $0\leq i \leq n-1$ and $x \in \mathbb{F}_{q^n}$. From Examples~\ref{ex:sup} and~\ref{ex22}, notice that $\displaystyle{\Psi^{-1}(N) = \frac{1}{\sigma_1 (\alpha)}\Psi^{-1}(\sigma_2 (M))}$, where $\sigma_i (M)$ is the a matrix obtained applying the automorphism $\sigma_i$ in each entry of $M$. This is a particular case of the next result.

\begin{proposition}\label{prop23}
Let $\sigma_j \in Aut\left(\mathbb{F}_{q^n}\right)$, for $1\leq j\leq n-1$. If $M \in \mathsf{Mat}_{m \times t}(\mathbb{F}_{q^n})$ is superregular, then $\sigma_j(M) \in \mathsf{Mat}_{m \times t} (\mathbb{F}_{q^n})$ is also superregular one.
\end{proposition}

\begin{proof}
It is straigthforward to see that 
\begin{equation*}
\begin{array}{cccc}
    \rho^{\sigma_j}:& \mathsf{Mat}_{m \times t}(\mathbb{F}_{q^n})  & \rightarrow &\mathsf{Mat}_{m \times t}(\mathbb{F}_{q^n})\\
    & M= \left(m_{i,j}\right) &\mapsto &  \rho^{\sigma_j}(M) = \sigma_j (M) = \left(\sigma_j(m_{i,j})\right) 
\end{array}    
\end{equation*}
is a ring isomorphism for each $1\leq j\leq n-1 $. Then, the result follows.
\end{proof}

Again, let $C$ be the companion matrix of an $n$-degree primitive polynomial in $\mathbb{F}_q[x]$.
Based on Corollary~\ref{co:Bi}, define $B_{\ell}=C^{k_{\ell}}$, with $\ell=1,\ldots, s$, where the exponents are pairwise different.
Then, the matrix (\ref{mat_prod_b2}) has the form:

\begin{equation}
M=\left[ 
\begin{array}{cccc}
a_{11}C^{k_1} & a_{12}C^{k_2}&\ldots &a_{1s}C^{k_s}\\
a_{21}C^{k_1}&a_{22}C^{k_2}&\ldots &a_{2s}C^{k_s}\\
\vdots &\vdots& &\vdots\\
a_{r1}C^{k_1}&a_{r2}C^{k_2}&\ldots&a_{rs}C^{k_s}
\end{array}
\right].
\end{equation}

Denote by $C^{t_{ij}}=a_{ij}C^{k_j}$,
for $i=1,\ldots, r$, $j=1, \ldots, s$, then the previous matrix is equal to:
\begin{equation*}
M=
\left[ 
\begin{array}{cccc}
C^{t_{11}} & C^{t_{12}}&\ldots &C^{t_{1s}}\\
C^{t_{21}}& C^{t_{22}}&\ldots &C^{t_{2s}}\\
\vdots &\vdots& &\vdots\\
C^{t_{s1}}&C^{t_
{s2}}&\ldots&C^{t_{ss}}
\end{array}
\right].
\end{equation*}

According to Corollary~\ref{co:Bi}, $M$ is an $n$-block superregular matrix. 
Therefore, given $\alpha\in \mathbb{F}_{q^n}$ a root of $p(x)$, the following matrix
\begin{equation*}
\Psi^{-1}(M)=
\left[ 
\begin{array}{cccc}
\alpha^{t_{11}} & \alpha^{t_{12}}&\ldots &\alpha^{t_{1s}}\\
\alpha^{t_{21}}& \alpha^{t_{22}}&\ldots &\alpha^{t_{2s}}\\
\vdots &\vdots& &\vdots\\
\alpha^{t_{s1}}&\alpha^{t_
{s2}}&\ldots&\alpha^{t_{ss}}
\end{array}
\right],
\end{equation*}
is superregular over $\mathbb{F}_{q^n}$. Notice that this method of constructing superregular matrices over extended alphabets is slightly more complex than the one presented before in this section. 

At next section, we propose a new ways to describe  several superregular matrices over a finite extension field, which might also be seen as either $n$-block superregular  matrices over $\mathbb{F}_q$, or superregular matrices over $GL_n (\mathbb{F}_q).$ 

\section{A variety of superregular matrices}\label{sec:var}

As expected, the larger the degree of extension field is, the more superregular matrices of a specific order we get. In this section, we present a method to obtain superregular matrices over a finite extension field from a given superregular matrix over the ground field.

First, we recall a minor well-known result
that will be useful in subsequent results. 

\begin{lemma}\label{lem:det}
    Let $A\in \mathsf{Mat}_n\left(\Fset_q\right)$ a square matrix whose rows are represented by $A_i$, $i=1, \ldots, n$.
    Assume that there exists
    $k\in [1,n]$, such that the $k-$th row is given by $A_k=\alpha X+\beta Y$, where $\alpha, \beta\in \Fset_q$ and $X,Y\in \Fset_q^n$, then
    $$
    \det(A) = 
\det
\left[
\begin{array}{c}
   A_1\\
   A_2\\
   \vdots\\
   A_{k-1}\\
   \alpha X+\beta Y\\
    A_{k+1}\\
    \vdots\\
    A_n\\
\end{array}
\right]=\alpha
\det
\left[
\begin{array}{c}
   A_1\\
   A_2\\
   \vdots\\
   A_{k-1}\\
    X \\
    A_{k+1}\\
    \vdots\\
    A_n\\
\end{array}
\right]+\beta
\det
\left[
\begin{array}{c}
   A_1\\
   A_2\\
   \vdots\\
   A_{k-1}\\
 Y\\
    A_{k+1}\\
    \vdots\\
    A_n\\
\end{array}
\right]
    $$
\end{lemma}

Now, we are ready to introduce one of the main results of this section.

\begin{theorem}\label{th40}
Let $M=\Psi^{-1}\left(\mathcal{M} \otimes C\right) \in \mathsf{Mat}_{m}\left(\mathbb{F}_{q^n}\right)$ be a superregular matrix, where $
\mathcal{M} \in \mathsf{Mat}_{m}\left(\mathbb{F}_{q}\right)$ is superregular and $C \in \mathsf{Mat}_{n}\left(\mathbb{F}_{q}\right)$ is the companion matrix of an $n-$degree primitive polynomial $p(x)\in \mathbb{F}_q [x]$. Given $\alpha \in \mathbb{F}_{q^n}$ a root of $p(x)$ and  $F\in \mathsf{Mat}_{m}\left(\mathbb{F}_{q^n}\right)$ such that 
\begin{equation*}
F= \left(\begin{array}{cccc}
0 & 0 & \ldots & 0  \\
\vdots & \vdots & \ldots & \vdots \\
f_1 & f_2 & \ldots & f_m\\
\vdots & \vdots & \ldots & \vdots \\
0 & 0 & \ldots & 0
\end{array}\right),   
\end{equation*}
where $f_j = \sum_{i=2}^n f_{j,i}\alpha^i \in \mathbb{F}_{q^n} $, for every $1\leq j\leq m$, the matrix $M + F$ is also superregular over $\mathbb{F}_{q^n}$. 
\end{theorem}

\begin{proof}
According to Theorems~\ref{thm:1} and~\ref{lemmasara}, let $M=\Psi^{-1}\left(\mathcal{M} \otimes C\right) \in \mathsf{Mat}_{m}\left(\mathbb{F}_{q^n}\right)$.
Since $M=[m_{i,j}]$ is superregular, then $\det(M) =\overline{m}\alpha^m = \det(\mathcal{M})\alpha^m \neq 0$; in particular, for any $j\times j-$submatrix $M_j$ of $M$, we have $ \det(M_j) =\overline{m_j}\alpha^j = \det(\mathcal{M}_j)\alpha^j \neq 0$, where $\mathcal{M}_j$ is the respective $j\times j-$submatrix of $\mathcal{M}$. 
We denote the respective submatrix obtained from $M$ after removing the $i-th$ row and $j-th$ column by $M_{i,j}$.

Without loss of generality, we may assume that the only nonzero row of $F$ is the first one, since for determinant computations the permutation of rows just change (or not) the sign of the former determinant. Let
\begin{equation*}
N=M +F =\left(\begin{array}{cccc}
m_{11} + \sum_{i=2} ^{n} f_{1,i} \alpha^{i}& m_{12} + \sum_{i=2} ^{n} f_{2,i} \alpha^{i}&\ldots  & m_{1m} + \sum_{i=2} ^{n} f_{m,i} \alpha^{i} \\
m_{21} &m_{22}&\ldots&m_{2m}\\
\vdots&\vdots&\ddots&\vdots\\
m_{m1} &m_{m2}&\ldots&m_{mm} 
\end{array}\right),  
\end{equation*}
where $f_{l,k} \in \mathbb{F}_q$ for $1\leq l\leq m$, $2\leq k\leq n$. According to Lemma~\ref{lem:det},
\begin{eqnarray}\label{naosei}
\det(N) &=& \det\left(M\right) + \det\left(\begin{array}{ccc}
    f_{1,2}\alpha^2 &\ldots & f_{m,2}\alpha^2 \\
    m_{2,1}&\ldots&m_{2,m}\\
    \vdots&\ddots&\vdots\\
    m_{m,1}&\ldots&m_{m,m}
\end{array}\right) 
+\cdots+ 
\det\left(\begin{array}{ccc}
    f_{1,n}\alpha^n &\ldots & f_{m,n}\alpha^n \\
    m_{2,1}&\ldots&m_{2,m}\\
    \vdots&\ddots&\vdots\\
    m_{m,1}&\ldots&m_{m,m}
\end{array}\right)\nonumber \\
&=&\overline{m}\alpha^m + \alpha^2 \left(\sum_{i=1}^m (-1)^{1+i} \cdot f_{i,2} \cdot \det(M_{1,i})\right)
+\cdots + 
\alpha^n \left(\sum_{i=1}^m (-1)^{1+i} \cdot f_{i,n} \cdot \det(M_{1,i})\right)\nonumber \\
&=&\overline{m}\alpha^m + \alpha^{2} \left(\sum_{i=1}^m (-1)^{1+i} \cdot f_{i,2} \cdot \overline{m_{1,i}}\cdot \alpha^{m-1}\right)
+\cdots+
\alpha^{n} \left(\sum_{i=1}^m (-1)^{1+i} \cdot f_{i,n} \cdot \overline{m_{1,i}}\cdot \alpha^{m-1}\right)\nonumber \\
&=&\overline{m}\alpha^m + \alpha^{m+1} \left(\sum_{i=1}^m (-1)^{1+i} \cdot f_{i,2} \cdot \overline{m_{1,i}}\right)+\cdots+ \alpha^{n+m-1} \left(\sum_{i=1}^m (-1)^{1+i} \cdot f_{i,n} \cdot \overline{m_{1,i}}\right)\nonumber\\
&=&\alpha^m \left[\overline{m} +  \left(\sum_{i=1}^m (-1)^{1+i} \cdot f_{i,2} \cdot \overline{m_{1,i}}\right)\alpha +\cdots+ \left(\sum_{i=1}^m (-1)^{1+i} \cdot f_{i,n} \cdot \overline{m_{1,i}}\right)\alpha^{n -1} \right]\nonumber \\
&=&\alpha^m \left[\overline{m} +  F_1 \alpha +\cdots+ F_{n-1}\alpha^{n -1} \right],
\end{eqnarray}
where $\displaystyle{F_{j-1} =\sum_{i=1}^m (-1)^{1+i} \cdot f_{i,j} \cdot \overline{m_{1,i}}}$, for $j\in \{2,\ldots,n\}$.

We state that~\eqref{naosei} is non-zero. It is enough to check that $\overline{m} +  F_1 \alpha +\cdots+ F_{n-1}\alpha^{n -1}\neq 0$. Indeed, since the polynomial $F(x) = \overline{m} + F_1 x + \cdots + F_{n-1} x^{n-1} \in \mathbb{F}_q [x]$ is nonzero, once $\overline{m}=\det(M)\neq 0$, and its degree is less  than $n$, the degree  of $p(x)$, which, by hypothesis,   is the minimal polynomial of $\alpha$, then it follows that $F(\alpha)\neq 0$  and, consequently $\alpha^m \left[\overline{m} +  F_1 \alpha +\cdots+ F_{n-1}\alpha^{n -1} \right]\neq 0.$

Now, we   consider a random  minor of $N$, denoted by $N_j$, of size $j\times j$  .
Let us assume that one of the $j$ rows of $N_j$ is the first row of $N$.
Otherwise, $N_j$ is non-singular, since it is a minor of the superregular matrix~$M$.

Let $\phi,\psi \in S_m$ be permutations so that $\phi(i)\neq 1$, for all $i\in\{2,\ldots,j\}$. Define
\begin{equation}
N_j = \left(\begin{array}{cccc}
    m_{1,\psi(1)} + \sum_{i=2} ^{n} f_{\psi(1),i} \alpha^{i}& m_{1,\psi(2)} + \sum_{i=2} ^{n} f_{\psi(2),i} \alpha^{i}&\ldots  & m_{1,\psi(j)} + \sum_{i=2} ^{n} f_{\psi(j),i} \alpha^{i} \\
m_{\phi(2),\psi(1)} &m_{\phi(2),\psi(2)}&\ldots&m_{\phi(2),\psi(j)}\\
\vdots&\vdots&\ddots&\vdots\\
m_{\phi(j),\psi(1)}&m_{\phi(j),\psi(2)}&\ldots&m_{\phi(j),\psi(j)}\\
\end{array}\right).    \nonumber 
\end{equation}

Without loss of generality, we denote by $N_{i,k}$ the respective submatrix of $N_j$ obtained after removing the $i-th$ row and $k-th$ column. Hence, 
\begin{eqnarray}\label{naosei2}
\det(N_j) &=& \det\left(M_j\right) \nonumber \\
&+& \det\left(\begin{array}{ccc}
    f_{\psi(1),2}\alpha^2 & \ldots & f_{\psi(j),2}\alpha^2 \\
    m_{\phi(2),\psi(1)}&\ldots&m_{\phi(2),\psi(j)}\\
    \vdots&\ddots&\vdots\\
m_{\phi(j),\psi(1)}&\ldots&m_{\phi(j),\psi(j)}
\end{array}\right)+\cdots+ \det\left(\begin{array}{ccc}
    f_{\psi(1),n}\alpha^n &\ldots & f_{\psi(j),n}\alpha^n \\
    m_{\phi(2),\psi(1)}&\ldots&m_{\phi(2),\psi(j)}\\
    \vdots&\ddots&\vdots\\
m_{\phi(j),\psi(1)}&\ldots&m_{\phi(j),\psi(j)}
\end{array}\right)\nonumber \\
&=&\overline{n_j}\alpha^j\nonumber \\
&+& \alpha^2 \left(\sum_{i=1}^j (-1)^{1+i} \cdot f_{\psi(i),2} \cdot \det(N_{1,i})\right)+\cdots + \alpha^n \left(\sum_{i=1}^j (-1)^{1+i} \cdot f_{\psi(i),n} \cdot \det(N_{1,i})\right)\nonumber \\
&=&\overline{n_j}\alpha^j \nonumber \\ 
&+& \alpha^{2} \left(\sum_{i=1}^j (-1)^{1+i} \cdot f_{\psi(i),2} \cdot \overline{n_{1,i}}\cdot \alpha^{j-1}\right) +\cdots+\alpha^{n} \left(\sum_{i=1}^j (-1)^{1+i} \cdot f_{\psi(i),n} \cdot \overline{n_{1,i}}\cdot \alpha^{j-1}\right)\nonumber \\
&=&\overline{n_j}\alpha^j + \alpha^{j+1} \left(\sum_{i=1}^j (-1)^{1+i} \cdot f_{\psi(i),2} \cdot \overline{n_{1,i}}\right)+\cdots+ \alpha^{n+j -1} \left(\sum_{i=1}^j (-1)^{1+i} \cdot f_{\psi(i),n} \cdot \overline{n_{1,i}}\right)\nonumber\\
    &=&\alpha^j \left[\overline{n_j} +  \left(\sum_{i=1}^j (-1)^{1+i} \cdot f_{\psi(i),2} \cdot \overline{n_{1,i}}\right)\alpha +\cdots+ \left(\sum_{i=1}^j (-1)^{1+i} \cdot f_{\psi(i),n} \cdot \overline{n_{1,i}}\right)\alpha^{n -1} \right]\nonumber \\
&=&\alpha^j \left[\overline{n_j} +  G_1 \alpha +\cdots+ G_{n-1}\alpha^{n -1} \right].
\end{eqnarray}

By a similar argument, we have that~\eqref{naosei2} is nonzero, and the result follows.
\end{proof}

\begin{example}\label{ex:tabelas}
Consider the superregular matrix in $\mathbb{F}_{13}$ given by:
$$
A=
\left[
\begin{array}{cccc}
   6 &  9  & 2&   6\\
   4  & 3 &  8&   1\\
   3  & 3&   4&   9\\
   3  & 9 &  9 &  5 \\
\end{array}
\right].
$$

Let $p(x)=x^3 +11x +6 \in \mathbb{F}_{13} [x]$ be a primitive polynomial, $\alpha\in \mathbb{F}_{13 ^3}$  one of its (primitive) roots, and $M=\Psi^{-1}(A\otimes C)$, where $C$ is the companion matrix of $p(x)$. 
According to Theorem~\ref{th40}, given $f  = f_2 \alpha^2 + f_3 \alpha^3$, $g = g_2 \alpha^2 + g_3 \alpha^3$, $h = h_2 \alpha^2 + h_3 \alpha^3$, and $i=i_2 \alpha^2 + i_3 \alpha^3 \in \mathbb{F}_{13^3}$, we can construct the following matrix
\begin{equation*}
 N=\left[
\begin{array}{cccc}
   6\alpha + f &  9\alpha +g  & 2\alpha + h &   6\alpha + i\\
   4\alpha  & 3\alpha &  8\alpha&   \alpha\\
   3\alpha  & 3\alpha&   4\alpha&   9\alpha\\
   3\alpha  & 9\alpha &  9\alpha &  5\alpha \\
\end{array}
\right]\in\mathsf{Mat}_{4}\left(\mathbb{F}_{13^3}\right).   
\end{equation*}

Notice that $\det(N)=\alpha^4 [(11 f_3 + 3 g_3 + 11 h_3+2 i_3)\alpha^2 +(11 f_2+3 g_2+11 h_2+2 i_2)\alpha +10]\neq 0$, and all the minors, described at Tables~\ref{table4} and~\ref{table3}, are different from zero as well.
Therefore, $N$ is a superregular matrix over $\mathbb{F}_{13^3}$.
\end{example}

Next, we introduce a second result that provides several novel superregular matrices over finite extension fields. Once again, our   construction depends solely on superregular matrices over the ground field, easily attainable, for instance, through classical constructions such as Vandermonde or Cauchy.

\begin{theorem}
Let $M=\Psi^{-1}\left(\mathcal{M} \otimes C\right)=[m_{i,j}]_m \in \mathsf{Mat}_{m}\left(\mathbb{F}_{q^n}\right)$ be a superregular matrix, where $
\mathcal{M} \in \mathsf{Mat}_{m}\left(\mathbb{F}_{q}\right)$ is superregular and $C \in \mathsf{Mat}_{n}\left(\mathbb{F}_{q}\right)$ the companion of a $n-$degree primitive polynomial $p(x)\in \mathbb{F}_q [x]$, respectively. Let $\alpha \in \mathbb{F}_{q^n}$ be a root of $p(x)$ and $\omega \in S_m$
so that $N_{\omega} \in \mathsf{Mat}_{m}\left(\mathbb{F}_{q^n}\right)$ is the matrix obtained permutating the rows of 
\begin{eqnarray*}
N&=& \left(\begin{array}{cccc}
 m_{1,1} + n_{1,1}\alpha^t & m_{1,2} + n_{1,2}\alpha^t & \ldots & m_{1,m}+ n_{1,m}\alpha^t  \\
 m_{2,1} + n_{2,1}\alpha^t & m_{2,2} + n_{2,2}\alpha^t & \ldots & m_{2,m}+ n_{2,m}\alpha^t  \\
\vdots & \vdots & \ddots & \vdots \\
 m_{j,1} + n_{j,1}\alpha^t & m_{j,2} + n_{j,2}\alpha^t & \ldots & m_{j,m}+ n_{j,m}\alpha^t  \\
 m_{j+1,1} & m_{j+1,2} & \ldots & m_{j+1,m} \\
\vdots & \vdots & \ddots & \vdots \\
m_{m,1} & m_{m,2} & \ldots & m_{m,m} \\
\end{array}\right) = \left(\begin{array}{c}
     \mathbf{m_1 + n_1\alpha^t}  \\
     \mathbf{m_2 + n_2\alpha^t}  \\ 
     \vdots\\
     \mathbf{m_j + n_j\alpha^t}  \\
     \mathbf{m_{j+1} }  \\ 
     \vdots\\
     \mathbf{m_{m} }  \\ 
\end{array}\right)\\
&=&\left(\begin{array}{c}
     (\mathbf{m+n})_{1,\ldots,j}  \\
     (\mathbf{m})_{j+1,\ldots,m}
\end{array}\right) 
\end{eqnarray*}
by $\omega.$ If $t>1$, $j>1$ and $j(t-1)< n$, then the matrix $N$ is  superregular over $\mathbb{F}_{q^n}$. 
\end{theorem}
\begin{proof}
Without loss of generality, we prove this theorem  for $N$, once the determinant of a row-permutation matrix obtained from $N$ is $det\left(N_{\omega}\right) = \pm \det(N)$. So,

\begin{eqnarray}\label{polinalpha}
\det(N)&=& \det\left(\begin{array}{c}
 \mathbf{m_{1}} \\
 \mathbf{m_{2}} \\
 \vdots\\
 \mathbf{m_{j}} \\
 \mathbf{m_{j+1}} \\
 \vdots\\
 \mathbf{m_{m}} \\
\end{array}\right) + \det\left(\begin{array}{c}
 \mathbf{n_{1}\alpha^t} \\
 \mathbf{m_{2}} \\
\vdots\\
 \mathbf{m_{j}} \\
 \mathbf{m_{j+1}} \\
 \vdots\\
 \mathbf{m_{m}} \\
\end{array}\right)+ \det\left(\begin{array}{c}
 \mathbf{m_{1}} \\
 \mathbf{n_{2}\alpha^t} \\
\vdots\\
 \mathbf{m_{j}} \\
 \mathbf{m_{j+1}} \\
 \vdots\\
 \mathbf{m_{m}} \\
\end{array}\right) +\det\left(\begin{array}{c}
 \mathbf{n_{1}\alpha^t}\\
 \mathbf{n_{2}\alpha^t} \\
\vdots\\
 \mathbf{m_{j}} \\
 \mathbf{m_{j+1}} \\
 \vdots\\
 \mathbf{m_{m}} \\
\end{array}\right) +\cdots \nonumber \\
&+&\det\left(\begin{array}{c}
 \mathbf{n_{1}\alpha^t}\\
 \mathbf{n_{2}\alpha^t} \\
\vdots\\
 \mathbf{n_{j}\alpha^t} \\
 \mathbf{m_{j+1}} \\
 \vdots\\
 \mathbf{m_{m}} \\
\end{array}\right),
\end{eqnarray}
and set $$N_{i_1 , i_2 , \ldots, i_k} = \left\{\begin{array}{cc}
   \mathbf{n_i}\alpha^t  ,&\mbox{if }i\in \{i_1 , i_2 , \ldots, i_k\}, \\
    \mathbf{m_i}\alpha^t ,&\mbox{otherwise, }
\end{array}\right.$$ where $\left\{i_1 , i_2 , \ldots, i_k\right\}\subset\{1,2,\ldots,j\}$. Thus,
\begin{eqnarray*}
\det(N)&=&\overline{m}\alpha^m  + \det(N_1)\alpha^t + \det(N_2)\alpha^t + \det(N_{1,2})\alpha^{2t} + \cdots + \det(N_{1,2,\ldots,j})\alpha^{jt}\nonumber \\
&=& \overline{m}\alpha^m  + \overline{n_1}\alpha^{m-1} \alpha^t + \overline{n_2}\alpha^{m-1} \alpha^t + \overline{n_{1,2}}\alpha^{m-2} \alpha^{2t} + \cdots + \overline{n_{1,2,\ldots,j}} \alpha^{m - j}\alpha^{jt} \nonumber\\
&=&\overline{m}\alpha^m  + \overline{n_1}\alpha^{m} \alpha^{t-1} + \overline{n_2}\alpha^{m} \alpha^{t-1} + \overline{n_{1,2}}\alpha^{m} \alpha^{2t-2} + \cdots + \overline{n_{1,2,\ldots,j}} \alpha^{m}\alpha^{j(t - 1)}\nonumber \\
&=&\alpha^m \left[\overline{m} + \left(\overline{n_1}+\cdots + \overline{n_j}  \right)\alpha^{t-1} +\cdots+ \left(\overline{n_{1,2,\ldots,j-1}} + \cdots +\overline{n_{2,\ldots,j}}\right)\alpha^{(j-1)(t-1)}\right. \nonumber \\ 
&+&\left.\overline{n_{1,2,\ldots,j}} \alpha^{j(t - 1)} \right].
\end{eqnarray*}   

We notice that the polynomial
\begin{equation*}
F(x)= \overline{m} + \left(\overline{n_1}+\cdots + \overline{n_j}  \right)x^{t-1} +\cdots+ \left(\overline{n_{1,2,\ldots,j-1}} + \cdots +\overline{n_{2,\ldots,j}}\right)x^{(j-1)(t-1)} 
 + \overline{n_{1,2,\ldots,j}} x^{j(t - 1)}   
\end{equation*}
is non-zero, since $\det(M)=\overline{m}\neq 0$ once $M$ is superregular. In addition, since the degree of $F(x)$ is $j(t-1)<n$, and the degree of the minimal polynomial of $\alpha$ is $n$, then $F(\alpha)\neq 0$. Thus, the element displayed in~\eqref{polinalpha} is non-zero, namely, $\det(N)\neq 0.$

Now, let $N_{\ell}$ be a $ {\ell}\times {\ell}$-minor of $N$, where $N_{\ell}$ composed by $r\geq 0$ rows of $(\mathbf{m+n})_{1,\ldots,j}$ and $s\geq 0$ rows of $(\mathbf{m})_{j+1,\ldots,m}$ such that ${\ell}=r+s$. Notice that, when $r=0$, then $N_{\ell}$ is a submatrix of
$M$, so it is non-singular. Thus, we can assume that $r\geq 1.$ So, in order to ensure the randomness of $N_{\ell}$ as a $\ell \times \ell$-submatrix of $N$, let $\psi \in S_j$, $\varphi \in S_{m - j}$ be permutations acting on the rows of $(\mathbf{m+n})_{1,\ldots,j}$ and $(\mathbf{m})_{j+1,\ldots,m}$, respectively,  and $\phi \in S_{\ell}$ is a permutation acting on columns of $N$.

For $\left\{i_1 , i_2 , \ldots, i_r \right\} \subset \{1 , 2, \ldots, j\}$, $\left\{j_1 , j_2 , \ldots, j_s \right\} \subset \{1 , 2, \ldots, m\}$ and $\left\{t_1 , t_2 , \ldots, t_{\ell} \right\} \subset \{1 , 2, \ldots, m\}$, a random $ {\ell}\times {\ell}$-minor of $N$ may be represented as
\begin{eqnarray*}
N_{\ell} &=& \left(\begin{array}{ccc}
m_{\psi(i_1),\phi(t_1)} + n_{\psi(i_1),\phi(t_1)}\alpha^t & \ldots& m_{\psi(i_1),\phi(t_{\ell})} + n_{\psi(i_1),\phi(t_{\ell})}\alpha^t  \\
m_{\psi(i_2),\phi(t_1)} + n_{\psi(i_2),\phi(t_1)}\alpha^t & \ldots& m_{\psi(i_2),\phi(t_{\ell})} + n_{\psi(i_2),\phi(t_{\ell})}\alpha^t  \\
\vdots & \ddots & \vdots \\
m_{\psi(i_r),\phi(t_1)} + n_{\psi(i_r),\phi(t_1)}\alpha^t &\ldots & m_{\psi(i_r),\phi(t_{\ell})} + n_{\psi(i_r),\phi(t_{\ell})} \alpha^t  \\
m_{\varphi(j_1),\phi(t_1)} & \ldots & m_{\varphi(j_1),\phi(t_{\ell})} \\
\vdots &  \ddots & \vdots \\
m_{\varphi(j_s),\phi(t_1)} &\ldots& m_{\varphi(j_s),\phi(t_{\ell})}
\end{array}\right),
\end{eqnarray*}
where 
$$M_{\ell}= \left(\begin{array}{c}
 \mathbf{m_{\psi(i_1)}} \\
  \vdots\\
 \mathbf{m_{\psi(i_r)}} \\
 \mathbf{m_{\varphi(j_1)}} \\
 \vdots\\
 \mathbf{m_{\varphi(j_s)}}\\
\end{array}\right)$$
and $\det(M_{\ell}) = \overline{m_{\ell}}\alpha^{\ell}\neq 0$, once $M_{\ell}$ is a submatrix of $M$, which is superregular.

Given a decomposition of the determinant of $N_{\ell}$ similar that one provided in~\eqref{polinalpha}, we have
\begin{eqnarray*}
\det(N_{\ell})&=&\overline{m_{\ell}}\alpha^{\ell}  + \det(N_{\psi(i_1)})\alpha^t + \det(N_{\psi(i_2)})\alpha^t + \det(N_{\psi(i_1), \psi(i_2)})\alpha^{2t} + \cdots +\nonumber \\
&+&\det(N_{\psi(i_1),\psi(i_2),\ldots,\psi(i_r) })\alpha^{rt}\nonumber \\
&=& \overline{m_{\ell}}\alpha^{\ell}  + \overline{n_{\psi(i_1)}}\alpha^{\ell-1} \alpha^t + \overline{n_{\psi(i_2)}}\alpha^{\ell-1} \alpha^t + \overline{n_{{\psi(i_1)},{\psi(i_2)}}}\alpha^{\ell-2} \alpha^{2t} + \cdots + \overline{n_{{\psi(i_1)},{\psi(i_2)},\ldots,{\psi(i_r)}}} \alpha^{\ell - r}\alpha^{rt} \nonumber\\
&=&\overline{m_{\ell}}\alpha^{\ell}  + \overline{n_{\psi(i_1)}}\alpha^{\ell + t-1} + \overline{n_{\psi(i_2)}}\alpha^{\ell + t -1}  + \overline{n_{{\psi(i_1)},{\psi(i_2)}}}\alpha^{\ell + 2t-2}  + \cdots + \overline{n_{{\psi(i_1)},{\psi(i_2)},\ldots,{\psi(i_r)}}} \alpha^{\ell + rt - r}
\nonumber \\
&=&\alpha^{\ell} \left[\overline{m}_{\ell} + \left(\overline{n_{\psi(i_1)}}+\cdots + \overline{n_{\psi(i_r)}}  \right)\alpha^{t-1} +\cdots+ \left(\overline{n_{\psi(i_1),\ldots,\psi(i_{r-1})}} + \cdots +\overline{n_{\psi(i_2),\ldots,\psi(i_{r})}}\right)\alpha^{(r-1)(t-1)}\right. \nonumber \\ 
&+&\left.\overline{n_{\psi(i_1),\ldots,\psi(i_r)}} \alpha^{r(t - 1)} \right].
\end{eqnarray*}   

As $r \leq  j$, then $r(t-1)\leq j(t-1)<n$,  and the polynomial
\begin{eqnarray*}
F_{\ell}(x)&=& \overline{m}_{\ell} + \left(\overline{n_{\psi(i_1)}}+\cdots + \overline{n_{\psi(i_r)}}  \right)\alpha^{t-1} +\cdots+ \left(\overline{n_{\psi(i_1),\ldots,\psi(i_{r-1})}} + \cdots +\overline{n_{\psi(i_2),\ldots,\psi(i_{r})}}\right)\alpha^{(r-1)(t-1)}\nonumber \\
&+& \overline{n_{\psi(i_1),\ldots,\psi(i_r)}} \alpha^{r(t - 1)}   
\end{eqnarray*}
is non-zero, 
since the degree of $F(x)$ is $j(t-1)<n$, and the degree of the minimal polynomial of $\alpha$ is $n$.
As a consequence,  $F(\alpha)\neq 0$ and thus, the determinant displayed in~\eqref{polinalpha} is non-zero, namely, $\det(N_{\ell})\neq 0.$
\end{proof}

  \begin{example} \label{ex:tabelas2}
Based on Example~\ref{ex:sup}, let $
\Psi^{-1}(M)=\left[\begin{array}{ccc}
    \alpha & 2\alpha & 2\alpha \\
    2\alpha & \alpha & 3\alpha \\
    3\alpha & 2\alpha & 4\alpha 
\end{array}\right]    \in \mathsf{Mat}_{3}(\mathbb{F}_{5^3})$ be a superregular matrix over $\mathbb{F}_{5^3}$, where $\alpha$ is a root of the primitive polynomial $p(x)=x^3 +3x +3\in \mathbb{F}_5[x]$.

Let $N=\left[\begin{array}{ccc}
    \alpha + n_{11}\alpha^2  & 2\alpha + n_{12}\alpha^2 & 2\alpha + n_{13}\alpha^2 \\
    2\alpha+ n_{21}\alpha^2 & \alpha + n_{22}\alpha^2& 3\alpha + n_{23}\alpha^2 \\
    3\alpha & 2\alpha & 4\alpha 
\end{array}\right]    \in \mathsf{Mat}_{3}(\mathbb{F}_{5^3}),$
where $n_{i,j}\in \mathbb{F}_5$, for $i=\{1,2\}$ and $j=\{1,2,3\}$. Notice that none of the entries of $N$ are zero. In addition,
\begin{eqnarray*}
\det(N) &=& \alpha^3\left[\mathbf{2} + (3n_{11}+n_{12}+n_{13}+n_{21}+3n_{22}+4n_{23})\alpha \right.\nonumber \\
&+& \left. (4n_{11}n_{22}+3n_{11}n_{23}+n_{12}n_{21}+3n_{12}n_{23}+2n_{13}n_{21}+2n_{13}n_{22})\alpha^2  \right]\neq 0    
\end{eqnarray*}

Finally,  all $2\times 2-$minors   are non-zero as shown in Table~\ref{table51}.
Therefore, $N$ is a superregular matrix over $\mathbb{F}_{5^3}$.
\end{example}

\section{Conclusions}\label{sec:conc}
In this paper we have introduced a construction of block superregular matrices based on the Kronecker product which can be reinterpreted via ring isomorphism as superregular matrices over finite extension fields. We have also provided two explicit methods to construct block superregular matrices which depend on from an unique superregular matrix over ground field. Considering large matrices, such a methods provide several superregular matrices over finite extension fields.

As future work, we would like to  explore the potential of sparse block superregular matrices due to the significant applications of sparse matrices in cryptography.
These matrices play a vital role in lattice-based cryptography, ensuring secure protocols, and are crucial in error-correcting codes for maintaining data transmission security. Moreover, their effectiveness in homomorphic cryptography enables the execution of operations on encrypted data without the need for decryption. Finally, we believe that some results presented in this paper may be used in order to construct lightweight matrices, which are also highly demanded in cryptography.
 
\section*{Acknowledgements}
The work of the second author was   supported by Conselho Nacional de Desenvolvimento Científico e Tecnológico (CNPq) with number of process  {405842/2023-6} and FAPESP with process  2024/04160-7.
 \bibliographystyle{elsarticle-num}
 \bibliography{biblio}
 
\newpage
 
 \appendix
 \section{}
 \begin{table}[h!]
\centering
\caption{All $3\times 3-$ minors of $N$ in Example~\ref{ex:tabelas}}
\label{table4}
\begin{tabular}{|c|c|}
\hline
\backslashbox{Rows}{Cols}  & $\{1,2,3\}$                                                                \\ \hline
$\{1,2,3\}$                                  & $\alpha^3[\textbf{6} +(f_2 + 8g_2 + 3h_2)\alpha + (f_3 + 8 g_3 + 3 h_3)\alpha^2] $       \\ \hline
$\{1,2,4\}$                                  & $\alpha^3[\textbf{1} +(7 f_2 + g_2+h_2)\alpha + (7 f_3 + g_3 + h_3)\alpha^2 ]$           \\ \hline
$\{1,3,4\}$                                  & $\alpha^3[\textbf{3} + (4f_2 + 11 g_2 + 5 h_2)\alpha + (4 f_3+11 g_3+5 h_3)\alpha^2 ]$   \\ \hline
\backslashbox{Rows}{Cols} & $\{1,2,4\}$                                                                \\ \hline
$\{1,2,3\}$                                  & $\alpha^3[\textbf{8} +(11 f_2 + 6 g_2 + 3 i_2)\alpha + (11f_3 + 6g_3 + 3i_3)\alpha^2 ]$  \\ \hline
$\{1,2,4\}$                                  & $\alpha^3[\textbf{6} +(6 f_2 + 9 g_2 + i_2)\alpha + (6f_3 + 9 g_3 + i_3)\alpha^2]$       \\ \hline
$\{1,3,4\}$                                  & $\alpha^3[\textbf{2} +(12 f_2 +12 g_2 + 5 i_2)\alpha + (12 f_3+12 g_3+5 i_3)\alpha^2 ]$  \\ \hline
\backslashbox{Rows}{Cols}  & $\{1,3,4\}$                                                                \\ \hline
$\{1,2,3\}$                                  & $\alpha^3[\textbf{8} +(3 f_2 + 6 h_2 + 5 i_2)\alpha + (3 f_3 + 6 h_3 + 5 i_3)\alpha^2 ]$ \\ \hline
$\{1,2,4\}$                                  & $\alpha^3[\textbf{1} +(7 f_2+g_2+h_2)\alpha + (7 f_3+g_3+h_3)\alpha^2 ]$                \\ \hline
$\{1,3,4\}$                                  & $\alpha^3[\textbf{8} +(4 f_2+12 h_2+2 i_2)\alpha +(4 f_3+12 h_3+2 i_3)\alpha^2 ]$       \\ \hline
\backslashbox{Rows}{Cols}  & $\{2,3,4\}$                                                                \\ \hline
$\{1,2,3\}$                                  & $\alpha^3 [\textbf{11} +(3 g_2+2 h_2+i_2)\alpha +(3 g_3 + 2 h_3 + i_3)\alpha^2 ]$       \\ \hline
$\{1,2,4\}$                                  & $\alpha^3[\textbf{1} + (7 f_2+g_2+h_2)\alpha + (7 f_3 + g_3+h_3)\alpha^2 ]$              \\ \hline
$\{1,3,4\}$                                  & $\alpha^3 [\textbf{10}+( 4 g_2+h_2+4 i_2)\alpha + (4 g_3+h_3+4 i_3)\alpha^2 ]$          \\ \hline
\end{tabular}
\end{table}
\begin{sidewaystable}
\centering
 \caption{All $2\times 2-$ minors of $N$  in Example~\ref{ex:tabelas}}
 \label{table3}
\begin{tabular}{|c|c|c}
\hline
\backslashbox{Rows}{Cols}          & $\{1,2\}$                                                 & \multicolumn{1}{c|}{$\{1,3\}$}                                                                                                                                      \\ \hline
$\{1,2\}$ & $\alpha^2[\textbf{8} +(3f_2+9g_2)\alpha + (3f_3+9g_3)\alpha^2 +(3f_2+9g_2)\alpha ]$        & \multicolumn{1}{c|}{$\alpha^2[\textbf{1} +(8f_2+9h_2)\alpha+
 (8f_3+9h_3)\alpha^2 ]$}                                                                                              \\ \hline
$\{1,3\}$ & $ \alpha^2 [\textbf{4}+ +(3f_2 + 10g_2)\alpha + (3f_3 + 10g_3) \alpha^2]$  & \multicolumn{1}{c|}{$\alpha^2[\textbf{5} +(4f_2 + 10h_2)\alpha + (4f_3 + 10h_3)\alpha^2 ]$} \\ \hline
$\{1,4\}$ & $\alpha^2 [\textbf{1} +(9f_2 + 10g_2)\alpha+ (9f_3 + 10g_3)\alpha^2  ]$ & \multicolumn{1}{c|}{$\alpha^2[\textbf{9} +(9f_2 + 10h_2)\alpha + (9f_3 + 10h_3)\alpha^2 ]$}                                                                                       \\ \hline
\backslashbox{Rows}{Cols}          & $\{1,4\}$                                                 & \multicolumn{1}{c|}{$\{2,3\}$}                                                                                                                                      \\ \hline
$\{1,2\}$ &  $\alpha^2 [\textbf{8} + (f_2 + 9 i_2) \alpha + (f_3 + 9i_3)\alpha^2  ]$                                                         &  \multicolumn{1}{c|}{$\alpha^2[\textbf{1} +(8 g_2 + 10 h_2) \alpha + (8 g_3+10 h_3) \alpha^2 ]$}                                                                                                                                               \\ \hline
$\{1,3\}$ &  $\alpha^2 [\textbf{10} + (9f_2+10 i_2) \alpha+ (9 f_3 + 10 i_3)\alpha^2 ]$                                                            & \multicolumn{1}{c|}{$\alpha^2[\textbf{4} +(4 g_2+10 h_2) \alpha + (4 g_3 + 10 h_3)\alpha^2 ]$}                                                                                                                                               \\ \hline
$\{1,4\}$ &   $\alpha^2 [\textbf{12} +(5 f_2 + 10 i_2)\alpha + (5 f_3 + 10 i_3) \alpha^2 ] $                                                       & \multicolumn{1}{c|}{$\alpha^2[\textbf{11} +(9 g_2 + 4 h_2) \alpha+ (9 g_3 + 4 h_3) \alpha^2 ]$}                                                                                                                                               \\ \hline
\backslashbox{Rows}{Cols}          & $\{3,4\}$                                                 & \multicolumn{1}{l}{}                                                                                                                                                \\ \cline{1-2}
$\{1,2\}$ &           $\alpha^2[ \textbf{6}+ (h_2 + 5i_2) \alpha  + (h_3 + 5 i_3) \alpha^2 ]$                                                & \multicolumn{1}{l}{}                                                                                                                                                \\ \cline{1-2}
$\{1,3\}$ &       $\alpha^2 [\textbf{7} +(9 h_2 + 9 i_2)\alpha + (9 h_3 + 9 i_3)\alpha^2]$                                                     & \multicolumn{1}{l}{}                                                                                                                                                \\ \cline{1-2}
$\{1,4\}$ &      $\alpha^2[\textbf{8}+ (5 h_2 + 4 i_2) \alpha +(5 h_3 + 4 i_3) \alpha^2 ]$                                                     & \multicolumn{1}{l}{} \\ \cline{1-2}
\end{tabular}
\end{sidewaystable}

\begin{table}[h!]
\centering
\caption{All $2\times 2-$ minors of $N$ in Example~\ref{ex:tabelas2}\label{table51}}
\begin{tabular}{|c|c|}
\hline
\backslashbox{Rows}{Cols}  & $\{1,2\}$                                                                \\ \hline
$\{1,2\}$                                  &  $\alpha^2 \left[\mathbf{2} + (n_{11} + 3n_{12} + 3n_{21} + n_{22})\alpha + (n_{11} n_{22} + 4n_{12} n_{21})\alpha^2 \right]$       \\ \hline
$\{1,3\}$                                  &      $\alpha^2\left[\mathbf{1} +(2n_{11} + 2n_{12})\alpha \right]$      \\ \hline
$\{2,3\}$                                  &   $\alpha^2 \left[\mathbf{1} + (2n_{21}+2n_{22})\alpha \right] $  \\ \hline
\backslashbox{Rows}{Cols} & $\{1,3\}$                                                                \\ \hline
$\{1,2\}$                                  & $ \alpha^2\left[\mathbf{4} + (3n_{11}+3n_{13}+3n_{21}+n_{23})\alpha+ (n_{11}n_{23}+4n_{13}n_{21})\alpha^2 \right]$  \\ \hline
$\{1,3\}$                                  &  $\alpha^2\left[\mathbf{3} + (4n_{11} + 2 n_{13})\alpha\right]$       \\ \hline
$\{2,3\}$                                  & $\alpha^2\left[\mathbf{4} +  (4n_{21}+2n_{23})\alpha\right]$   \\ \hline
\backslashbox{Rows}{Cols}  & $\{2,3\}$                                                                \\ \hline
$\{1,2\}$                                  &  $\alpha^2\left[\mathbf{4} + (3n_{12}+4n_{13}+3n_{22}+2n_{23})\alpha + (n_{12}n_{23}+4n_{13}n_{22})\alpha^2 \right]$ \\ \hline
$\{1,3\}$                                  &  $\alpha^2\left[\mathbf{4} +  (4n_{12} + 3n_{13})\alpha \right]$                \\ \hline
$\{2,3\}$                                  &   $\alpha^2\left[\mathbf{3} + (4n_{22}+3n_{23})\alpha \right]$     \\ \hline
\end{tabular}
\end{table}

\end{document}